\documentclass[12pt]{article}

\usepackage{graphicx, epsfig}
\usepackage{amsmath,amssymb}
\usepackage[authoryear]{natbib}
\usepackage[latin1]{inputenc}

\begin{document}

\title{Small-sample likelihood inference in extreme-value regression models}
\author{Silvia L.P. Ferrari
%$^{\ast}$\thanks{$^\ast$Corresponding author. Email: silviaferrari.usp@gmail.com}
\qquad  Eliane C. Pinheiro\\
{\small {\em Departamento de Estat\' \i stica, Universidade de S\~ao Paulo,
%\\Rua do Mat\~ao, 1010, S\~ao Paulo/SP, 05508-090, 
Brazil}}
}

\maketitle
\begin{abstract}

We deal with a general class of extreme-value regression models introduced by
Barreto-Souza and Vasconcelos (2011). Our goal is to derive an adjusted likelihood ratio statistic that 
is approximately distributed as $\chi^2$ with a high degree of accuracy. Although the 
adjusted statistic requires more computational effort than its unadjusted counterpart, it is shown that the adjustment term has a simple compact form that 
can be easily implemented in standard statistical software. 
Further, we compare the finite sample 
performance of the three classical tests (likelihood ratio, Wald, and score), the gradient test that has been recently proposed by \cite{TERRELL}, and the adjusted likelihood ratio test 
obtained in this paper. Our simulations favor the latter. Applications of our results are presented.


\end{abstract}

\section{Introduction}

The extreme-value distributions are frequently used to model extreme events, such as extreme floods
and wind speed, and in survival or reliability analysis to model the logarithm of lifetime data.
The literature on statistics of extremes has grown fast due to the increasing 
interest in statistically modeling extreme values (minimum or maximum) in a wide range of areas, such as
climatology, hydrology, reliability, finance, insurance, and environmental sciences. A classical
reference on the subject is the book by \cite{GUMBEL}. More recent references include 
\cite{KOTZNADA}, \cite{COLES}, and \cite{CASTILLO}, among others.

In this paper, we deal with a general class of extreme-value regression models introduced by
\cite{BARRETO}. The authors presented large-sample inference on the parameters, and also 
considered the issue of correcting the bias of the maximum likelihood estimators in small samples. 
The study of small-sample inference in extreme-value models is relevant since
the amount of extreme data available for analysis may be small in practical applications and, as 
\cite{MUELLERRUFIBACH} pointed out, there are very few articles that focus on small sample problems in extreme value theory.
Here, we focus on statistical tests in the general class of extreme-value regression
models proposed by \cite{BARRETO} when the sample size is small.
Specifically, our goal is to derive Skovgaard's adjusted likelihood ratio statistics in this 
class of models. We show that the adjustment term 
has a simple compact form that can be easily implemented from standard statistical software.
The adjusted statistic is approximately distributed as $\chi^2$ with a high degree of
accuracy. 
Further, we compare the finite sample performance of the three classical tests ( likelihood ratio, Wald, and score), the gradient test that has been recently proposed 
by \cite{TERRELL}, and the adjusted likelihood ratio test obtained in this paper. 

Let $y_1, \ldots, y_n$ be independent random variables, where 
each $y_t$, $t=1,\ldots,n$, has an extreme-value distribution with parameters $\mu_t$ and
$\phi_t$ and density function
\begin{equation}
\label{extreme-value}
f(y; \mu_t, \phi_t) = \frac{1}{\phi_t}\exp\left(-\frac{y-\mu_t}{\phi_t}\right)\exp\left\{-\exp\left(-\frac{y-\mu_t}{\phi_t}\right)\right\} ,\quad y \in I\!\! R,
\end{equation}
where $\mu \in I\!\! R$ and $\phi > 0$ are the location and dispersion parameters, respectively.
The mean and the variance of $y_t$ are
${\rm E}(y_t) = \mu_t+{\mathcal E}\phi_t$
and 
${\rm var}(y_t) = {\phi_t^2\pi^2}/{6},$
respectively, where ${\mathcal E}$ is the Euler constant; ${\mathcal E}\approx 0.5772$. 
If $y$ has an extreme-value distribution with parameters $\mu$ and $\phi$, we write $y\sim EV_{max}(\mu,\phi)$.
This distribution is also called Gumbel or type I extreme-value distribution. Here, we call it the maximum 
extreme-value distribution to contrast it with the minimum extreme-value distribution (see Section 3 below), 
which is also called Gumbel or type I extreme-value in the statistical literature.

The maximum extreme-value regression model with dispersion covariates is defined by (\ref{extreme-value}) and by two systematic components given by
\begin{equation}\label{linkmu}
g(\mu_t) = \eta_t = \eta(x_{t},\beta)
\end{equation}
and
\begin{equation}\label{linkphi}
h(\phi_t) =  \delta_t = \delta(z_{t},\gamma),
\end{equation}
where $\beta = (\beta_1, \ldots, \beta_k)^{\!\top}$ and $\gamma = (\gamma_1, \ldots, \gamma_m)^{\!\top}$ are vectors of  unknown regression parameters ($\beta \in I\!\! R^k$ and $\gamma \in I\!\! R^m$, $k + m < n$) and 
$x_{t}$ and $z_{t}$ are observations on covariates.
Here, $\eta(\cdot,\cdot)$ and $\delta(\cdot,\cdot)$ are continuously twice differentiable 
(possibly nonlinear) functions in the second argument.
Finally, $g(\cdot)$ and $h(\cdot)$ are known strictly monotonic and twice differentiable link functions that map $I\!\! R$ and $I\!\! R^+$ 
into $I\!\! R$, respectively.
Let $X$ be the derivative matrix of $\eta=(\eta_1,\ldots,\eta_n)^{\top}$ with respect to $\beta^\top$. 
Analogously, let $Z$ be the derivative matrix of $\delta=(\delta_1,\ldots,\delta_n)^{\top}$ with respect to $\gamma^\top$. 
It is assumed that ${\rm rank}(X)=k$ and ${\rm rank}(Z)=m$ for all $\beta$ and $\gamma$.

The paper unfolds as follows. 
In Section \ref{main}, we derive Skovgaard's adjusted likelihood ratio statistic for testing hypothesis on the parameters of the model definide by (\ref{extreme-value})-(\ref{linkphi}). 
In Section \ref{min-extreme-value}, we extend our results to a general class
of minimum extreme-value regression models. 
Section \ref{MonteCarlo} is devoted to a simulation study to compare the performance
of the three classical tests, the gradient test, and the adjusted likelihood ratio test. Our
simulation results clearly favor the adjusted test proposed in this paper. 
In Section \ref{App}, we illustre the use of our 
results in two real data sets. Section \ref{conclusion} closes the paper with a discussion. Technical details are left for an appendix.

\section{Main results}\label{main}

Let $\ell(\theta)$ be the log-likelihood function of the model defined by
(\ref{extreme-value})-(\ref{linkphi}) 
given the vector of observations $y=(y_1,\ldots,y_n)$. We have
$$
\ell(\theta) = \sum_{t=1}^n \ell_t(\mu_t, \phi_t),
$$
where 
$$
\ell_t(\mu_t, \phi_t) = -\log (\phi_t) - \frac{y_t-\mu_t}{\phi_t} 
- \exp\left( - \frac{y_t-\mu_t}{\phi_t}\right),
$$
with $\mu_t$ and $\phi_t$ defined so that (\ref{linkmu}) and (\ref{linkphi}) hold. In matrix notation, the log-likelihood
function can be written as
\begin{equation}\label{loglikvet} 
\ell(\theta) = [-\mathnormal {l}^{\!\top}-{\mathfrak z}^{\!\top}-{\breve {\mathfrak z}}^{\!\top}]{\mathbf 1},
\end{equation}
where
${\mathnormal l} = (\ln \phi_1,\ldots,\ln \phi_n)^{\!\top}$, 
${\mathfrak z}=({\mathfrak z} _1,\ldots,{\mathfrak z}_n)^{\!\top}$, and
${\breve {\mathfrak z}}=(\exp(-{\mathfrak z} _1),\ldots,\exp(-{\mathfrak z}_n))^{\!\top}$, with ${\mathfrak z}_t=(y_t-\mu _t)/\phi _t$,
and 
${\mathbf 1}$ is the $n$-dimensional column vector of ones.
The score function, obtained by differentiating the log-likelihood function with respect  to the unknown parameters, is denoted by 
$U \equiv (U_\beta(\beta,\gamma)^\top,U_\gamma(\beta,\gamma)^\top)^\top$, with
$I$ and $J$ denoting the expected and observed information matrices. We have
\begin{equation}\label{scorebeta}
U_{\beta}(\beta, \gamma) = X^{\!\top} \Phi^{-1} T ({\mathbf 1}-{\breve {\mathfrak z}}),
\end{equation}
\begin{equation}\label{scoregamma}
U_\gamma (\beta,\gamma )= Z^{\!\top} \Phi^{-1} H ({\mathfrak z}-{\mathcal Z}{\breve {\mathfrak z}}-{\mathbf 1}),
\end{equation}
$$
J =\left[
\begin{array}{c c}
J_{\beta \beta } & J_{\beta \gamma }\\
J_{\gamma \beta } & J_{\gamma \gamma }
\end{array}
\right],
\qquad
I=\left [
\begin{array}{c c}
I_{\beta\beta}&I_{\beta\gamma}\\
I_{\gamma\beta}&I_{\gamma\gamma}
\end{array}
\right],
$$
%\hspace{0.1cm} \rm{with} \hspace{0.1cm}
%\left\{
with
\begin{eqnarray*}
J_{\beta \beta }&=&X^{\!\top}\Phi ^{-1}T\Bigl({\breve {\mathcal Z}}\Phi ^{-1} + ({\mathcal I}-{\breve {\mathcal Z}})ST\Bigr)TX-\Bigl[{\mathbf 1}^{\!\top}\bigl({\mathcal I}-{\breve {\mathcal Z}}\bigr)T\Phi ^{-1}\Bigr]\Bigl[\dot X\Bigr], \\
J_{\beta \gamma }&=&J_{\gamma \beta }^{\!\top}=X^{\!\top}\Phi^{-1}T\bigl({\mathcal I}-{\breve {\mathcal Z}}+{\mathcal Z}{\breve {\mathcal Z}}\bigr)H\Phi^{-1}Z, \\
J_{\gamma \gamma }&=&
Z^{\!\top}\Phi^{-1}H\Bigl((-{\mathcal I}+2{\mathcal Z}-2{\mathcal Z}{\breve {\mathcal Z}}+{\mathcal Z}^2{\breve {\mathcal Z}})\Phi^{-1} \\
& &+(-{\mathcal I}+{\mathcal Z}-{\mathcal Z}{\breve {\mathcal Z}})QH\Bigr)HZ
+\Bigl[{\mathbf 1}^{\!\top}\bigl({\mathcal I}-{\mathcal Z}+{\mathcal Z}{\breve {\mathcal Z}}\bigr)H\Phi^{-1}\Bigr]\Bigl[\dot Z\Bigr], \\
I_{\beta \beta }&=&X^{\!\top}\Phi ^{-2}T^2X, \quad 
I_{\beta \gamma }=I_{\gamma \beta }^{\!\top}=({\mathcal E}-1)X^{\!\top}\Phi^{-1}TH\Phi^{-1}Z, 
%\quad {\rm and} 
\\
I_{\gamma \gamma }&=&\bigl(1 + \Gamma ^{(2)}(2)\bigr)Z^{\!\top}\Phi^{-1}H^2\Phi^{-1}Z,
\end{eqnarray*}
where
%${\mathcal Z}\;=\;{\rm diag}\{{\mathfrak z} _1,\ldots,{\mathfrak z}_n\}$,
%${\breve {\mathcal Z}}\;=\;{\rm diag}\{\exp(-{\mathfrak z}_1),\ldots,\exp (-{\mathfrak z} _n)\}$,
%$\Phi\;=\;{\rm diag}\{\phi_1,\ldots,\phi_n\}$,
%$T\;=\;{\rm diag}\{ 1/g'(\mu_1), \ldots, 1/g'(\mu_n)\}$,
%$H\;=\;{\rm diag}\{ 1/h'(\phi_1), \ldots, 1/h'(\phi_n)\}$,
%$S\;=\;{\rm diag}\{g''(\mu _1),\ldots,g''(\mu _n) \}$,
%$Q\;=\;{\rm diag}\{h''(\phi _1),\ldots,h''(\phi _n) \}$,
${\mathcal I}$ is an $n \times n$ identity matrix,
${\mathcal Z}={\rm diag}\{{\mathfrak z} _1,\ldots,{\mathfrak z}_n\}$,
${\breve {\mathcal Z}}={\rm diag}\{\exp(-{\mathfrak z}_1),\ldots,\exp (-{\mathfrak z} _n)\}$,
$\Phi={\rm diag}\{\phi_1,\ldots,\phi_n\}$,
$T={\rm diag}\{ 1/g'(\mu_1), \ldots, 1/g'(\mu_n)\}$,
$H={\rm diag}\{ 1/h'(\phi_1), \ldots, 1/h'(\phi_n)\}$,
$S={\rm diag}\{g''(\mu _1),\ldots,g''(\mu _n) \}$,
$Q={\rm diag}\{h''(\phi _1),\ldots,h''(\phi _n) \}$,
$\dot X=\partial ^2 \eta / \partial \beta \partial \beta ^{\top}$ and
$\dot Z=\partial ^2 \delta / \partial \gamma \partial \gamma ^{\top}$
are $n \times k \times k$ and $n \times m \times m$ arrays, respectively, and
$\Gamma(\cdot)$ denotes the gamma function, with
$\Gamma ^{(1)}(\cdot)$ and $\Gamma ^{(2)}(\cdot)$ being its first and second derivatives,
respectively. The other quantities are as before.

Let $\theta=(\beta^\top,\gamma^\top)^\top$ be the unknown parameter vector
that indexes the extreme-value regression model (\ref{extreme-value})-(\ref{linkphi}). 
In what follows, $\nu=(\nu_1,\ldots,\nu_r)^\top$ represents the parameter of interest  
and $\psi=(\psi_1,\ldots,\psi_s)^\top$ is the nuisance parameter; note that $r+s=k+m$. 
We consider likelihood-based tests of the null hypothesis  ${\mathcal H}_0: \nu=\nu_0$, 
where $\nu_0$ is a fixed  $r$-vector. Clearly, such tests may be inverted to give confidence 
sets for $\nu$.
Further, let $J_{\psi \psi}$ denote the $s \times s$
observed information matrix corresponding to $\psi$. Similarly, 
$A_{\psi \psi}$ denotes a matrix formed from the $(r+s) \times (r+s)$ matrix 
$A$ by dropping the rows and columns that correspond to the interest parameter.
Additionally, hat and tilde indicate evaluation at the unrestricted 
($\widehat\theta$) and at the restricted ($\widetilde\theta$) maximum likelihood 
estimator of $\theta$ under ${\mathcal H}_0$, respectively. For instance,
$\widehat{I}=I(\widehat{\theta})$, $\widetilde{I}=I(\widetilde{\theta})$, and
$\widehat{J}=J(\widehat{\theta})$. 

\cite{SKOVGAARD} derived  an adjusted likelihood ratio statistic given by
\begin{equation}\label{w*}
w^{*}=w - 2 \; {\rm log}\; {\zeta},
\end{equation}
where 
$
w =  2 (\ell(\widehat\theta)-\ell(\widetilde\theta))
$
is the likelihood ratio statistic,
$$
{\zeta}=
\frac{\{\mid\!\!\widetilde{I}\!\!\mid \;
\mid\!\!\widehat{I}\!\!\mid \;
\mid\!\!\widetilde{J}_{\psi\psi}\!\!\mid\}^{1/2}}
{\mid\!\!\overline{\Upsilon}\!\!\mid \; \mid\!\!\{\widetilde{I}
\overline{\Upsilon}^{-1} \widehat{J} \widehat{I}^{-1}
\overline{\Upsilon}\}_{\psi\psi}\!\!\mid^{1/2}}
\frac{ \{\widetilde{U}^{\top} \overline{\Upsilon}^{-1} \widehat{I}
\widehat{J}^{-1} \overline{\Upsilon} \widetilde{I}^{-1} \widetilde{U}\}^{r/2} }
{w^{r/2-1} \widetilde{U}^{\top} \overline{\Upsilon}^{-1} \overline{q}}, 
$$
and $\overline q$ and $\overline \Upsilon$ come, respectively, from  
\begin{equation}\label{q}
q={\rm E}_{\theta_1}[U(\theta_1) \ (\ell(\theta_1)-\ell(\theta))] 
\end{equation}
and
\begin{equation}\label{Upsilon}
\Upsilon={\rm E}_{\theta_1}[U(\theta_1) \ U^{\top}(\theta)]
\end{equation}
by inserting  $\widehat \theta$ for $\theta_1$ and $\widetilde \theta$
for $\theta$ after the expected values are computed. Note that
$\overline q$ is an $(r+s)$-vector and $\overline \Upsilon$ is an
$(r+s) \times (r+s)$ matrix.
Under ${\mathcal H}_0$, $w$ is distributed as $\chi^2_r$ with error of order $n^{-1}$ while
$w^{*}$ follows this distribution with high degree of accuracy \citep[p.~7]{SKOVGAARD}.
Simulation results in \cite{FERRARI2010} and \cite{FERRARI2008} suggest that tests that use $w^*$ are much less size distorted
than those that are based on $w$.

In order to obtain the adjusted likelihood ratio statistic (\ref{w*}) 
in the extreme-value regression model (\ref{extreme-value})-(\ref{linkphi}), one needs to obtain 
the score vector, the observed and expected information matrices, $J$ and $I$, respectively, the 
vector $q$, and the matrix $\Upsilon$. We obtained
\vspace{0.3cm}
$$%\label{qoverline}
\overline q = 
\left[
\begin{array}{c}
\widehat X^\top \widehat \Phi ^{-1} \widehat T  C ({\mathcal I} - M{\breve D}) {\mathbf 1} \\
\widehat Z^\top \widehat \Phi ^{-1} \widehat H \{ C({\mathcal E}{\mathcal I} + N {\breve D})-{\mathcal I} \}{\mathbf 1}
\end{array}
\right]
$$
and
$$%\label{Upsilonoverline}
\overline \Upsilon=
\left[
\begin{array}{c c}
\widehat X^\top \widehat \Phi ^{-1}\widehat T C M {\breve D} \widetilde T \widetilde \Phi^{-1} \widetilde X &

\widehat X^\top \widehat \Phi ^{-1} \widehat T  C\{ {\mathcal I}+{\breve D}(MD-M-CN)\}  \widetilde H \widetilde \Phi ^{-1} \widetilde Z \\

-\widehat Z^\top \widehat \Phi ^{-1} \widehat H C N {\breve D}  \widetilde T \widetilde \Phi ^{-1} \widetilde X  \quad&
\widehat Z^\top \widehat \Phi ^{-1} \widehat H  C \{{\mathcal E} {\mathcal I} + {\breve D}(N+CP-ND)\}  \widetilde H \widetilde \Phi ^{-1} \widetilde Z
\end{array}
\right],
$$
where
$C\;=\;{\rm diag}\{\phi_{11}/\phi_1,\ldots,\phi_{1n}/\phi_n\}$, 
$D\;=\;{\rm diag}\{(\mu_{11}-\mu_1)/\phi_1,\ldots,(\mu_{1n}-\mu_n)/\phi_n\}$,
${\breve D}\;=\;{\rm diag}\{\exp(-(\mu_{11}-\mu_1)/\phi_1),\ldots,\exp(-(\mu_{1n}-\mu_n)/\phi_n)\}$,
$M\;=\;{\rm diag}\{\Gamma (1+\phi_{11} / \phi_1),\ldots,\Gamma (1+\phi_{1n} / \phi_n)\}$,
$N\;=\;{\rm diag}\{\Gamma ^{(1)} (1+\phi_{11} / \phi_1),\ldots,\Gamma ^{(1)} (1+\phi_{1n} / \phi_n)\}$,
$P\;=\;{\rm diag}\{\Gamma ^{(2)} (1+\phi_{11} / \phi_1),\ldots,\\ \Gamma ^{(2)} (1+\phi_{1n} / \phi_n)\}$,
and the other quantities are as given above. Details of the derivations of $\overline q$ and $\overline \Upsilon$ are
given in the Appendix.%~\ref{ape:apendixqUpsilon}.

\section{Minimum extreme-value regression model}\label{min-extreme-value}

Let $y_1, \ldots, y_n$ be independent random variables, where 
each $y_t$, $t=1,\ldots,n$, has a minimum extreme-value distribution with parameters $\mu_t$ and
$\phi_t$ and density function
\begin{equation}
\label{eq:min-extreme-value}
f(y; \mu_t, \phi_t) = \frac{1}{\phi_t}\exp\left(\frac{y-\mu_t}{\phi_t}\right)\exp\left\{-\exp\left(\frac{y-\mu_t}{\phi_t}\right)\right\} ,\quad y \in I\!\! R,
\end{equation}
where $\mu \in I\!\! R$ and $\phi > 0$ are the location and dispersion parameters, respectively.
The mean and the variance of $y_t$ are
${\rm E}(y_t) = \mu_t-{\mathcal E}\phi_t$
and 
${\rm var}(y_t) = {\phi_t^2\pi^2}/{6},$
respectively. If $y$ has  a minimum extreme-value distribution we write $y\sim EV_{min}(\mu,\phi)$.
A useful property of the minimum extreme-value distribution is as follows:
\begin{eqnarray}\label{min-max}
y\sim EV_{min}(\mu,\phi) \Longrightarrow -y\sim EV_{max}(-\mu,\phi).
\end{eqnarray}
The minimum extreme-value regression model with dispersion covariates is defined by (\ref{eq:min-extreme-value}) and by the systematic components
(\ref{linkmu}) and (\ref{linkphi}). 

From (\ref{min-max}), it is easy to see that the minimum extreme-value regression model (\ref{eq:min-extreme-value}) with systematic components (\ref{linkmu}) and (\ref{linkphi}) is equivalent to the (maximum) extreme-value regression model (\ref{extreme-value}) for the response variables 
$v_1=-y_1,\ldots,v_n=-y_n$ with 
systematic components $g^*(\mu_t)=\mu_t^*=\eta^*_t=-g^{-1}(\eta(x_t,\beta))$ and $h^*(\phi_t)=h(\phi_t)=\delta(z_t,\gamma)$. Hence,   
inference for the minimum extreme-value regression model (\ref{eq:min-extreme-value}) with systematic components (\ref{linkmu}) and (\ref{linkphi}) 
can be performed from the results derived in Section \ref{main} by changing the signs of the observations on the response variable and using an identity link function for the location parameter with the modified predictor $\eta^*_t$.
As a result, the adjusted likelihood ratio statistic derived in Section \ref{main} can be easily computed for the minimum extreme-value regression model.

\section{Monte Carlo simulation results}\label{MonteCarlo}

We now present Monte Carlo simulation results on the small sample behaviour of the likelihood ratio
test ($w$), the Wald test ($W$), the score test ($S_R$), the gradient test ($S_T$),  and the adjusted 
likelihood ratio test ($w^*$). 
The Wald, score, and gradient statistics are given by
$W=(\widehat \nu - \nu_0)^\top \; (\widehat I^{\nu\nu})^{-1} \; (\widehat \nu - \nu_0),$
$S_R=\widetilde U_{\nu}^\top \; \widetilde I^{\nu\nu} \; \widetilde U_{\nu}$, and 
$S_T= \widetilde U_{\nu}^\top(\widehat \nu - \nu_0)$. 
Note that the gradient statistic is very simple to compute since it does not involve the information matrix, 
neither the observed one nor the expected one.

The maximum likelihood estimation of the unknown parameters was performed using the quasi-Newton BFGS nonlinear optimization 
algorithm with analytical derivatives developed by Broyden, Fletcher, Goldfarf \& Shanno (see, for instance, \cite{PRESS}) and 
implemented in the function MAXBFGS in the matrix language programming {\tt Ox} \citep{DOORNIK}.

All the size simulation results are based on 10,000 Monte Carlo replications and the nominal level of the tests are $\alpha = 10\%$, 5\%, and 1\%.
We also present power simulation results. Since the different tests display different sizes when a $\chi^2$ distribution is used, 
we simulated 500,000 samples to estimate the critical values of the tests that give exact size, i.e., size equal to the chosen nominal level.
Our power simulation results are obtained using exact critical values.

%\subsection[linear precisão constante]{Modelo de regressão valor extremo máximo linear com dispersão constante}
%\label{linear_homo}

We consider model (\ref{extreme-value}) with constant dispersion and
location parameter given by
$$
\mu_t = \beta_1 + \beta_2 x_{t2} + \beta_3 x_{t3}
                                    + \beta_4 x_{t4} + \beta_5 x_{t5},
$$
wich we refer as `model 1'.
Three null hypotheses are considered, ${\mathcal H}_0: \beta_2=0$ ($r=1$),
${\mathcal H}_0: \beta_2=\beta_3=0$ ($r=2$), and ${\mathcal H}_0: \beta_2=\beta_3=\beta_4=0$
($r=3$), and these are to be tested against a two-sided alternative.
For the first case, we set  $\beta_1=1$, $\beta_2=0$, $\beta_3=1$, $\beta_4=6$, and $\beta_5=-3$; for the second case,
$\beta_1=1$, $\beta_2=\beta_3=0$, $\beta_4=6$, and $\beta_5=-3$; and for the third case, 
$\beta_1=1$, $\beta_2=\beta_3=\beta_4=0$, and $\beta_5=-3$.
The value of $\phi$ was fixed at $\phi=0.1$.\footnote{For linear extreme-value regression models with constant dispersion, 
the null distributions of the five statistics do not depend on $\phi$. The proof is omitted to save space.}
%; see Appendix~\ref{ape:apendixDistribution}}. 
The covariate values were obtained as random draws from a ${\mathcal U}(-0.5, 0.5)$ distribution and the sample 
sizes are 15, 20, 30, and 40.

Table~\ref{tab:linear-homo_size} presents the null rejection rates of the five tests. 
It can be noticed that the likelihood ratio and Wald tests are markedly liberal in small samples.
For instance, for $n=15$, $r=1$, and $\alpha=5\%$, the null rejection rates of these tests are 11.6\% and 20.8\%,
respectively. 
The gradient test is liberal in many cases but not as much as the likelihood ratio and the Wald tests. 
The score test is less liberal and displays conservative behavior in some cases.
The adjusted likelihood ratio test is clearly the least size distorted. 
For the case mentioned above, the null rejection rates of the score, gradient, and adjusted likelihood ratio tests are 
5.8\%, 9.0\%, and 5.0\%, respectively.

\begin{table}[htp]
\centering
\renewcommand{\arraystretch}{1.1}
\caption{Null rejection rates (\%); model 1}\label{tab:linear-homo_size}
\footnotesize
\begin{tabular}
{c|c| c c c c c | c c c c c|  c c c c c}
\hline \hline
& & \multicolumn{5}{c|}{$\alpha = 10\%$}
& \multicolumn{5}{c|}{$\alpha = 5\%$}
& \multicolumn{5}{c}{$\alpha = 1\%$}
\\
\cline{3-17}
$r$ & $n$ & $w$ & $W$ & $S_R$ & $S_T$ & $w^*$ %& $w^{**}$ 
          & $w$ & $W$ & $S_R$ & $S_T$ & $w^*$ %& $w^{**}$ 
          & $w$ & $W$ & $S_R$ & $S_T$ & $w^*$ %& $w^{**}$
            \\
\hline            
 &15	&19.3	&28.2	&11.5	&16.9	&10.4	&11.6	&20.8	&5.8	&9.0	&5.0	&4.0	&11.1	&1.1	&1.8	&1.0 \\
1&20	&16.8	&22.1	&11.3	&15.0	&9.9	&10.0	&15.6	&5.4	&7.8	&4.8	&3.0	&7.0	&0.9	&1.5	&0.9 \\
 &30	&14.9	&18.9	&10.9	&13.8	&10.6	&8.5	&12.0	&5.3	&7.1	&5.4	&2.4	&4.6	&0.9	&1.4	&1.2 \\
 &40	&13.0	&15.7	&10.3	&12.1	&10.1	&7.0	&9.5	&4.9	&6.2	&5.1	&2.0	&3.5	&0.9	&1.3	&1.1 \\
\hline															
 &15	&22.9	&36.5	&12.7	&16.4	&7.7	&14.3	&28.7	&6.2	&7.9	&3.5	&4.8	&17.1	&1.1	&0.8	&0.6 \\
2&20	&18.8	&28.8	&11.6	&14.1	&9.2	&11.4	&20.9	&5.4	&7.2	&4.6	&3.8	&10.5	&0.9	&0.8	&0.8 \\
 &30	&16.0	&23.4	&10.8	&13.0	&10.2	&9.2	&15.7	&5.2	&6.6	&5.2	&2.3	&6.7	&0.8	&0.9	&0.9 \\
 &40	&13.9	&19.0	&10.2	&11.7	&9.8	&7.6	&12.0	&5.0	&6.1	&5.2	&1.9	&4.5	&0.7	&0.9	&0.9 \\
\hline															
 &15	&23.4	&42.9	&9.1	&12.6	&6.8	&14.9	&34.2	&3.4	&4.8	&3.0	&5.1	&21.2	&0.2	&0.2	&0.5 \\
3&20	&19.8	&33.8	&9.8	&12.3	&8.7	&12.1	&25.1	&4.4	&5.4	&4.4	&3.8	&14.0	&0.5	&0.4	&0.8 \\
 &30	&16.6	&26.4	&10.1	&12.1	&9.9	&9.5	&18.2	&4.7	&5.4	&4.9	&2.5	&8.6	&0.6	&0.6	&0.9 \\
 &40	&14.7	&21.4	&10.2	&11.3	&10.3	&8.3	&14.0	&4.6	&5.5	&5.2	&2.1	&5.8	&0.7	&0.7	&0.9 \\
\hline \hline
\end{tabular}
\end{table}

\begin{figure}[!h]
  \centering
  \includegraphics[height=150mm,width=165mm]{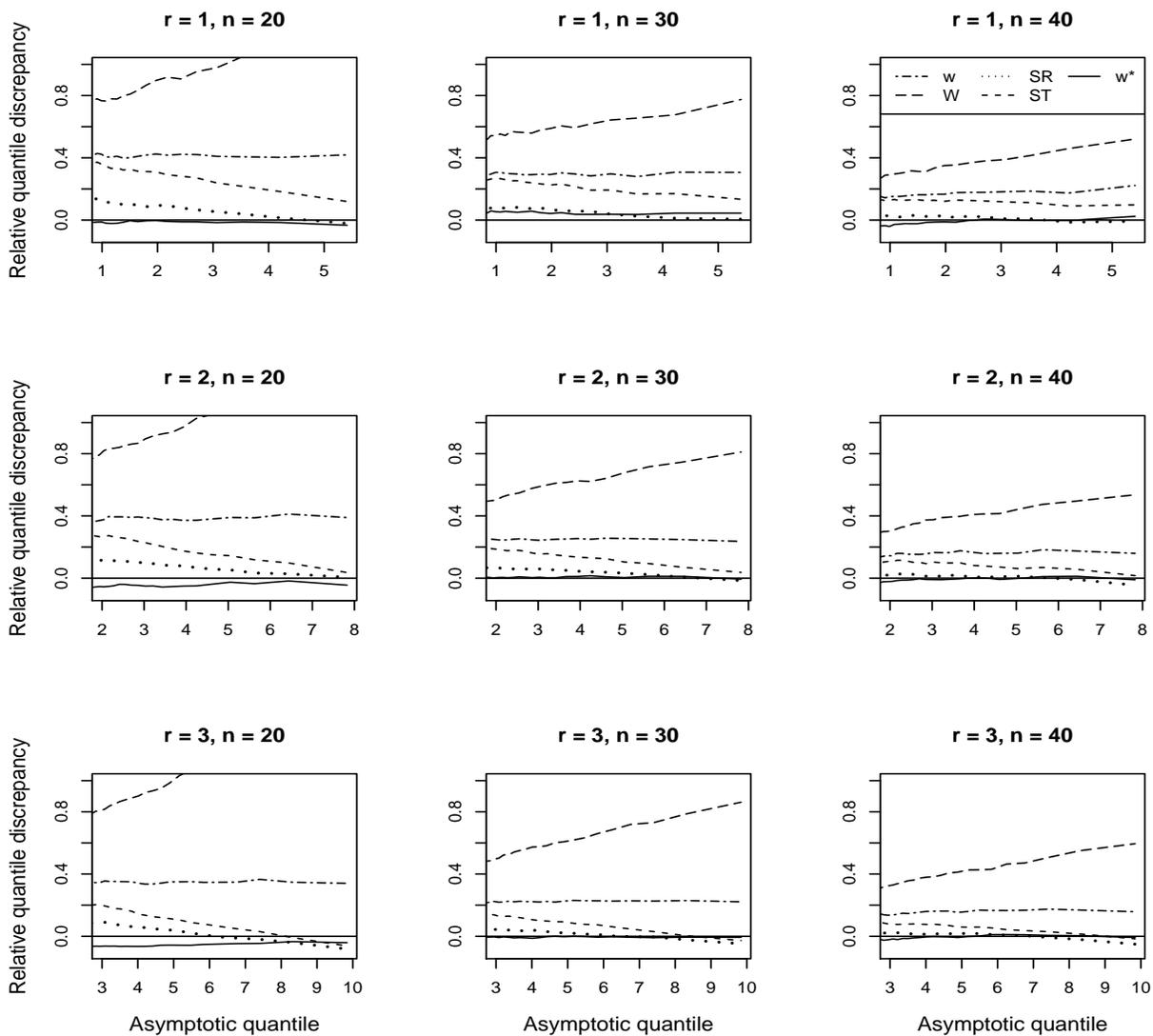}
  \caption{Relative quantile discrepancies, model 1.}
  \label{fig:linear-homo}
\end{figure}

Figure~\ref{fig:linear-homo} shows the plots of the relative quantile discrepancies versus corresponding asymptotic quantiles for $r=1,2,3$,
and $n=20, 30, 40$. Relative quantile discrepancy is defined as the difference between exact (estimated by simulation) and asymptotic quantiles divided by the latter. The closer to zero the relative quantile discrepancy, the better is the approximation of the exact null distribution of the test statistic by the limiting $\chi^2$ distribution. 
The plots confirm the tendency of the likelihood ratio and the Wald tests of rejecting the null hypothesis with higher frequency than expected based 
on the nominal level. It is clear that the distribution of the adjusted likelihood ratio statistic ($w^*$) closely agrees with the reference 
distribution. The effect of the proposed adjustment becomes evident. Notice that the best agreement between the exact distribution and its 
asymptotic counterpart is achieved by $w^*$ for all the sample sizes.

\begin{figure}[!h]
  \centering
  \includegraphics[height=77mm,width=127mm]{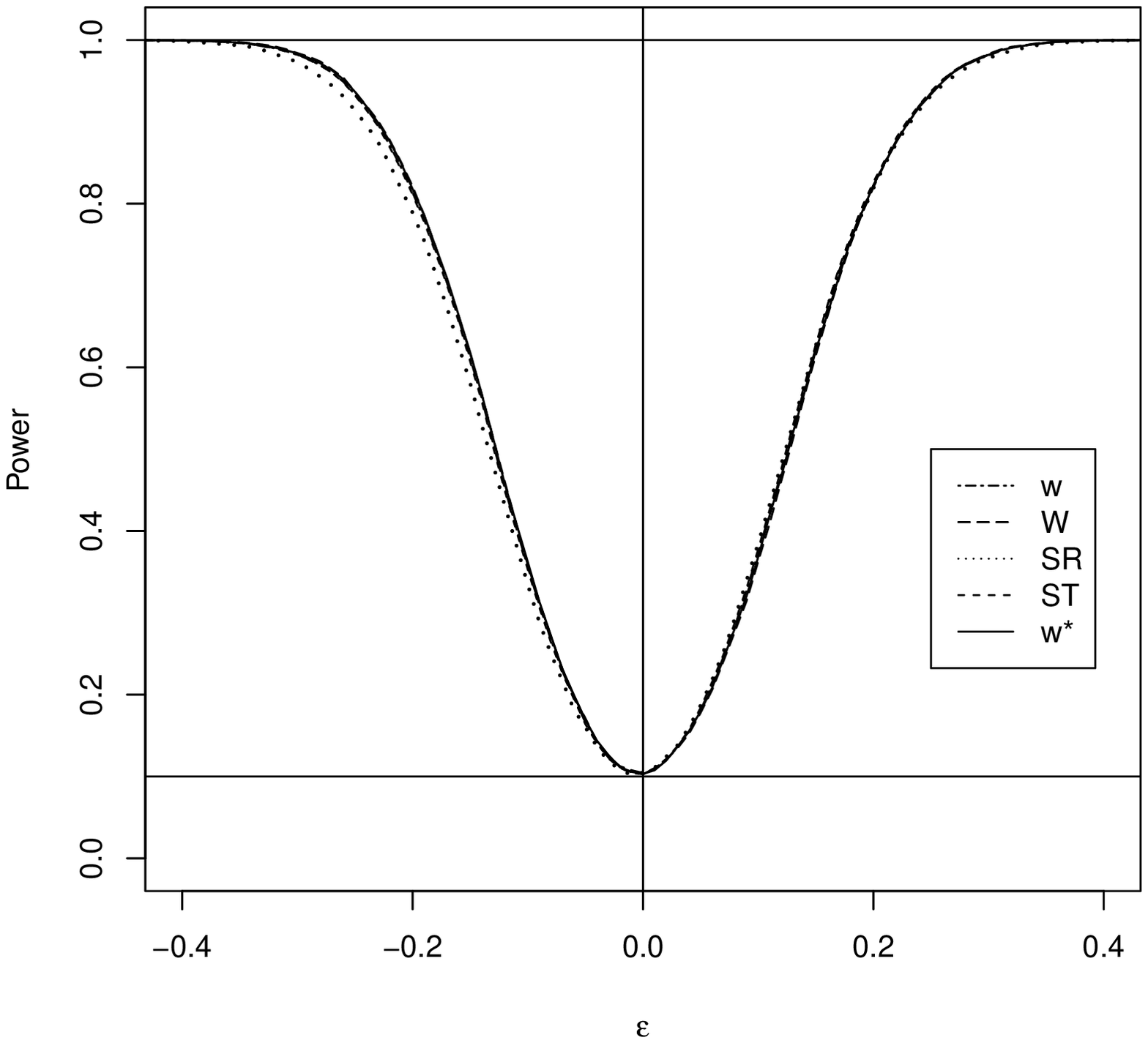}
  \caption{Power of the tests; model 1; $r=1$, $n=30$, $\alpha=10\%$.}
  \label{fig:linear-homo-r1-poder-500mil-n30}
\end{figure}

%%%%%%%%%%%%%%%%%%%%%%%%%%%%%%%%%%%%%%   poder
We now focus on the power comparisons of the five tests for $r=1$, $n=30$, and $\alpha=10\%$. 
The rejection rates were obtained under the alternative hypothesis 
${\mathcal H}_1: \beta_2=\epsilon$ for different values of $\epsilon$ through Monte Carlo simulation.
Figure~\ref{fig:linear-homo-r1-poder-500mil-n30} gives the plots of the power function of the tests. Visual inspection shows that 
the curves are practicaly coincident, i.e., the five tests display similar powers.

%{Modelo de regressão valor extremo linear com precisão variável}
%{Teste dos coeficientes da regressão do parâmetro de locação}

Now, consider model (\ref{extreme-value}) with systematic components for the location and scale parameters, respectively, given by 
$$
\mu_t = \beta_1 + \beta_2 x_{t2} + \beta_3 x_{t3} + \beta_4 x_{t4},
$$
and
$$
\ln (\phi_t) = \gamma_1 + \gamma_2 z_{t2} + \gamma_3 z_{t3} + \gamma_4 z_{t4},
$$
wich we refer to as `model 2'.

We consider three different null hypotheses, ${\mathcal H}_0: \beta_4=0$ ($r=1$),
${\mathcal H}_0: \beta_3=\beta_4=0$ ($r=2$), and ${\mathcal H}_0: \beta_2=\beta_3=\beta_4=0$ ($r=3$), and these are to be tested against two-sided alternatives.
The values for the $\beta$s are $\beta_1=1$, $\beta_2=1$, $\beta_3=6$, and $\beta_4=0$;
$\beta_1=1$, $\beta_2=1$, and $\beta_3=\beta_4=0$; and
$\beta_1=1$ and $\beta_2=\beta_3=\beta_4=0$, for the first, second, and third cases, respectively.
Further, we set $\gamma_1=\ln(0.1)\approx -2.30$, %-2.30259
$\gamma_2=-2$, $\gamma_3=-2$, and $\gamma_4=0.1$.
The covariate values were randomly drawn from a ${\mathcal U}(-0.5, 0.5)$ distribution and the sample sizes are 40, 50, 60, and 70.
%Table~\ref{tab:linear-hetero-H0beta_tamanho} \ref{tab:linear-hetero-H0beta_poder} e na
%Figura~\ref{fig:linear-hetero-H0beta}.

\begin{figure}%[!h]
  \centering
  \includegraphics[height=150mm,width=165mm]{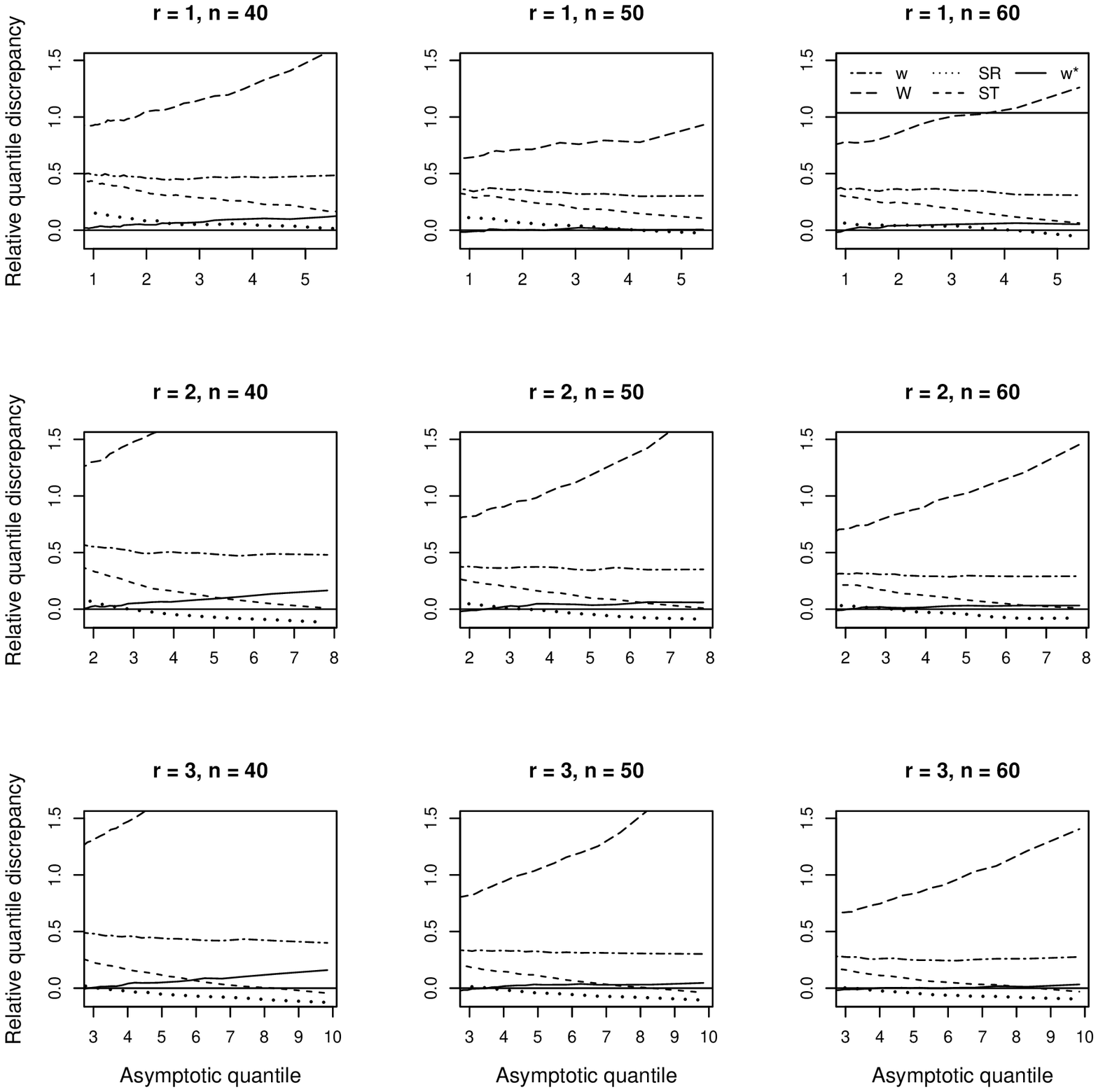}
  \caption{Relative quantile discrepancy; model 2.}
  \label{fig:linear-hetero-H0beta}
\end{figure}

Table~\ref{tab:linear-hetero-H0beta_size} gives the null rejection rates of the five tests and Figure~\ref{fig:linear-hetero-H0beta}  
shows the plots of relative quantile discrepancies. 
We note that the results for model 2 show similarity with those for model 1. 
The likelihood ratio and Wald tests are clearly oversized, i.e. its type I error probability is greater than the nominal level, 
and their distributions are much different from the asymptotic $\chi^2$ distribution if the sample is not large. 
Again, the best agreement between the true and asymptotic quantiles is reached by $w^*$, the adjusted likelihood ratio statistic
proposed in this paper. The score test presents good behavior but tends to be conservative when $r>1$. 
The gradient test is liberal in many cases, but high order quantiles of the gradient statistic are close to the asymptotic quantiles. 
We also performed power simulation comparisons among the five tests. Overall, the tests are equally powerful when true critical values are used.
\begin{table}%[htp]
\centering
\renewcommand{\arraystretch}{1.1}
\caption{Null rejection rates (\%); model 2}\label{tab:linear-hetero-H0beta_size}
\footnotesize
\begin{tabular}
{c|c| c c c c c | c c c c c|  c c c c c}
\hline \hline
& & \multicolumn{5}{c|}{$\alpha = 10\%$}
& \multicolumn{5}{c|}{$\alpha = 5\%$}
& \multicolumn{5}{c}{$\alpha = 1\%$}
\\
\cline{3-17}
$r$ & $n$ & $w$ & $W$ & $S_R$ & $S_T$ & $w^*$ %& $w^{**}$ 
          & $w$ & $W$ & $S_R$ & $S_T$ & $w^*$ %& $w^{**}$ 
          & $w$ & $W$ & $S_R$ & $S_T$ & $w^*$ %& $w^{**}$
            \\
\hline            
1	&	40	&	17.5	&	24.2	&	10.8	&	15.3	&	11.1	&	10.4	&	16.8	&	5.7	&	8.4	&	6.1	&	3.4	&	8.0	&	1.0	&	1.7	&	1.6	\\
	&	50	&	15.9	&	20.6	&	10.8	&	14.1	&	10.1	&	8.9	&	13.4	&	5.2	&	7.2	&	5.3	&	2.3	&	5.4	&	0.8	&	1.4	&	1.1	\\
	&	60	&	15.9	&	21.9	&	10.6	&	13.9	&	10.7	&	9.3	&	15.3	&	5.1	&	7.0	&	5.7	&	2.4	&	7.0	&	0.8	&	1.2	&	1.2	\\
	&	70	&	14.8	&	19.9	&	10.0	&	12.9	&	10.5	&	8.4	&	13.1	&	5.0	&	6.6	&	5.4	&	2.0	&	5.6	&	0.7	&	1.1	&	1.2	\\
\hline
2	&	40	&	21.7	&	36.9	&	8.3	&	13.7	&	11.8	&	13.7	&	28.6	&	3.7	&	6.4	&	6.6	&	4.5	&	16.6	&	0.6	&	0.9	&	2.0	\\
	&	50	&	18.6	&	29.4	&	9.0	&	13.5	&	10.9	&	11.0	&	21.4	&	3.9	&	6.3	&	5.7	&	3.2	&	11.1	&	0.5	&	0.9	&	1.4	\\
	&	60	&	17.0	&	26.9	&	9.0	&	12.8	&	10.6	&	9.7	&	19.6	&	3.9	&	5.9	&	5.4	&	2.9	&	9.9	&	0.5	&	0.9	&	1.3	\\
	&	70	&	15.8	&	24.5	&	8.8	&	12.4	&	10.3	&	9.0	&	17.1	&	4.0	&	5.6	&	5.2	&	2.3	&	8.0	&	0.6	&	0.9	&	1.2	\\
\hline
3	&	40	&	23.1	&	43.2	&	7.9	&	12.1	&	11.9	&	14.1	&	34.6	&	3.3	&	5.1	&	6.8	&	4.7	&	21.7	&	0.3	&	0.6	&	2.1	\\
	&	50	&	19.3	&	34.0	&	8.1	&	12.1	&	10.9	&	11.5	&	25.9	&	3.5	&	5.2	&	5.6	&	3.4	&	14.4	&	0.4	&	0.7	&	1.4	\\
	&	60	&	17.2	&	30.4	&	8.1	&	11.4	&	10.1	&	9.9	&	22.3	&	3.6	&	5.4	&	5.2	&	2.8	&	11.6	&	0.5	&	0.8	&	1.3	\\
	&	70	&	16.4	&	27.3	&	8.3	&	11.4	&	10.2	&	9.4	&	19.8	&	3.6	&	5.0	&	5.0	&	2.3	&	9.2	&	0.6	&	0.7	&	1.1	\\
\hline \hline
\end{tabular}
\end{table}

%{Modelo de regressão valor extremo não linear com precisão constante}

We now consider model (\ref{extreme-value}) with constant dispersion and a nonlinear specification for the location parameter
given by 
$$
\mu_t = \beta_0 + \beta_1 x_{t1} +x_{t2}^ {\beta_2},
$$
wich we refer to as `model 3'.

Here, $X$ is an $n\times 3$ matrix whose $t$-th row is $(1,x_{t1},\ln(x_{t2})x_{t2}^{\beta _2})$, $Z \! = \! {\mathcal I}$,
$T \! = \! {\mathcal I}$, $H\! = \! {\mathcal I}$, $S$ and $Q$ are matrices of zeros,
$\dot X= \partial ^2 \eta / \partial \beta \partial \beta ^{\top}$ is such that
$\partial ^2 \eta_t / \partial \beta_2 \partial \beta_2=\ln(x_{t2})^2x_{t2}^{\beta _2}$ for $t=1,\ldots,n$ and zero otherwise, and
$\dot Z = \partial ^2 \delta / \partial \gamma \partial \gamma ^{\top}$
is an $n \times m \times m$ array of zeros.
The bracket product of the $1 \times n$ vector $\Bigl[{\mathbf 1}^{\!\top}\bigl({\mathcal I}-{\breve {\mathcal Z}}\bigr)T\Phi ^{-1}\Bigr]$ and
the $n \times 3 \times 3$ array $\dot X$ is an $1 \times 3 \times 3$ array, i.e., a $3 \times 3$ matrix, whose $(i,j)$-th
element is 
$$
\sum_{t=1}^n\Biggl\{
\left(\frac{1}{\phi _t}-\frac{1}{\phi _t}\exp \left(-\frac{y_t-\mu _t}{\phi _t}\right)\right)\frac{1}{g'(\mu_t)}\ln(x_{t2})^2x_{t2}^{\beta _2}\Biggr\}
$$
if $(i,j)=(3,3)$ and zero otherwise.

\begin{figure}%[!h]
  \centering
  \includegraphics[height=50mm,width=165mm]{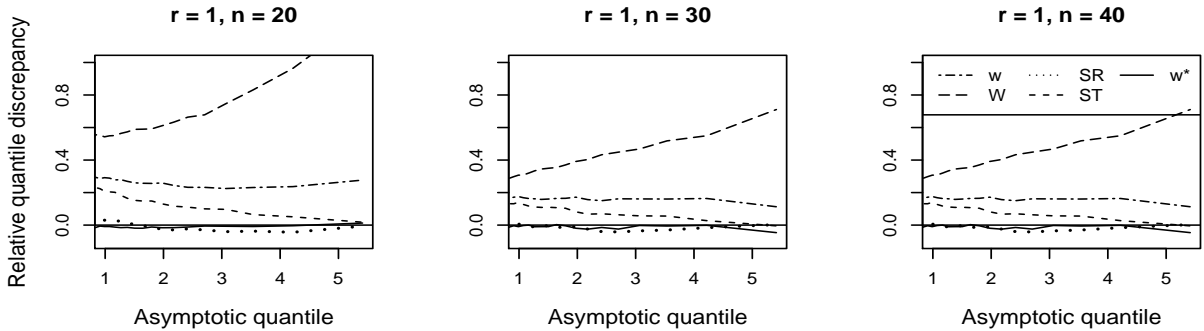}
  \caption{Relative quantile discrepancy, model 3.}
  \label{fig:nonlinear-homo}
\end{figure}

The null hypothesis under test is ${\mathcal H}_0: \beta_2=0$ ($r=1$).
We set $\phi={\rm e}^{0.1}\approx 1.1$, %\approx 1.105171 
$\beta_0=1$, $\beta_1=1$, and $\beta_2=0$, and 
the covariate values were drawn from a ${\mathcal U}(0,1)$ distribution.
The sample sizes are 15, 20, 30, and 40.
Table~\ref{tab:nonlinear-homo_size} and Figure~\ref{fig:nonlinear-homo} show our simulation results.
\begin{table}%[htp]
\centering
\renewcommand{\arraystretch}{1.1}
\caption{Null rejection rates (\%); model 3}\label{tab:nonlinear-homo_size}
\footnotesize
\begin{tabular}
{c| c c c c c | c c c c c|  c c c c c}
\hline \hline
& \multicolumn{5}{c|}{$\alpha = 10\%$}
& \multicolumn{5}{c|}{$\alpha = 5\%$}
& \multicolumn{5}{c}{$\alpha = 1\%$}
\\
\cline{2-16}
$n$ & $w$ & $W$ & $S_R$ & $S_T$ & $w^*$ %& $w^{**}$ 
    & $w$ & $W$ & $S_R$ & $S_T$ & $w^*$ %& $w^{**}$ 
    & $w$ & $W$ & $S_R$ & $S_T$ & $w^*$ %& $w^{**}$
\\
\hline
15	&	16.6	&	22.2	&	10.2	&	13.4	&	10.0	&	9.8	&	15.8	&	4.6	&	6.3	&	5.0	&	2.9	&	8.2	&	0.9	&	1.1	&	1.1	\\
20	&	14.1	&	19.2	&	9.5	&	11.9	&	9.9	&	7.7	&	12.7	&	4.5	&	5.7	&	4.8	&	2.2	&	5.8	&	1.1	&	1.0	&	1.1	\\
30	&	12.5	&	16.3	&	9.3	&	11.1	&	9.7	&	6.8	&	10.3	&	4.6	&	5.7	&	5.1	&	1.5	&	3.9	&	1.3	&	1.0	&	0.9	\\
40	&	12.2	&	15.1	&	9.8	&	11.2	&	10.3	&	6.7	&	9.3	&	4.9	&	5.6	&	5.2	&	1.5	&	3.4	&	1.4	&	1.1	&	1.0	\\
\hline \hline
\end{tabular}
%\end{center}
\end{table}

From Table~\ref{tab:nonlinear-homo_size}, we note that the tests that use  $w$ and $W$ are typically liberal while 
the null rejection of the other tests keeps their sizes closer to the nominal levels. The adjusted likelihood ratio test and 
the score test display better performance than the others.
For example, for $n=15$ and $\alpha =10\%$, the null rejection rates of the tests are $16.6\%$ (likelihood ratio), $22.2\%$ (Wald),
$13.4\%$ (gradient), $10.2\%$ (score), and $10.0\%$ (adjusted likelihood ratio).

Figure~\ref{fig:nonlinear-homo} shows that the reference distribution is not a good approximation for the null distribution of $w$ and $W$, 
but is close to the true null distribution of the score and the adjusted likelihood ratio statistics. 
The high order quantiles of the gradient statistic closely agree with the corresponding quantiles of the reference distribution.
The five tests have similar power performance (results not shown).

\section{Applications}\label{App}

In this section, we illustrate applications of our results in two real datasets. First, we deal with a dataset 
presented in Faivre \& Masle~(1988) and in Huet et al.~(2004). The aim is to study the growth of winter wheat, 
by measuring the differences in dry weights of wheat tillers and stems. The explanatory variable $x$, measured 
on a cumulative degree-days scale, is an integral in time of all temperatures at which the wheat is submitted 
that are above the smallest temperature at which wheat can develop. Temperatures are measured in degrees Celsius 
and time is measured in days, the initial time being determined by the physiological state of the wheat. 
Plants growing on $n = 18$ randomly chosen small areas of about $0.15$ $m^2$ are harvested
each week and the dry weights of the tillers for plants harvested from each area are measured in milligrams.
A detailed description of the data can be found in \citet[p. 61]{HUET}.

\cite{BARRETO} assumed that the dry weight of tillers ($y_1,\ldots,y_{18}$) is independent and follows a
nonlinear extreme-value regression model (\ref{extreme-value}) with 
$$
\mu_t=\beta_0+{\rm e}^{\beta_1+\beta_2x_t}
\ \ {\rm and} \ \ 
\ln \phi_t=\gamma_1x_t, \ \  t=1,\dots,18.
$$
They focused on the issue of correcting the bias of the maximum likelihood estimates (MLEs) of the parameters.
The authors constructed confidence intervals based on the asymptotic normality of the MLEs and of the 
bias-corrected MLEs. Their simulation study suggested that the asymptotic confidence intervals centered 
at the bias-corrected estimators produce coverage probability closer to the nominal confidence coefficient 
than those centered in the uncorrected MLEs. However, the choice of the estimator (uncorrected or corrected MLE) 
does not change the approximation error between the true coverage probabilty and the nominal confidence coefficient. 
In fact, the correction that they derived only guarantees that the bias of the corrected estimators are of 
order $O(n^{-2})$, but does not change the convergence rate of the distribution of the estimators to the normal 
distribution. 
%They obtained the following asymptotic 95\% confidence 
%intervals based on the bias-corrected MLEs: $CI(\beta_0,95\%)=(48,4 \; ; \; 111,0)$, $CI(\beta_1,95\%)=(-4,71\; ;\; -0,0673)$, $CI(\beta_2,95\%)=(0,0110\; ;\; 0,0169)$ and $CI(\gamma_0,95\%)=(0,00622\; ;\; 0,00738)$.
 
Here, we illustrate the use of the five test statistics, namely the likelihood ratio, Wald, score, gradient, 
and adjusted statistic derived in this paper, to obtain interval estimates for the parameters. By choosing 
the 5\% nominal level, the approximate confidence coefficient is 95\%. We emphasize that the confidence intervals 
obtained from the adjusted likelihood ratio statistic have coverage probabilities that are approximately equal to 
95\% with a high degree of accuracy. This is guaranteed by the theoretical results in \cite{SKOVGAARD} and 
is confirmed in our simulation study.

\begin{figure}[!h]
  \centering
  \includegraphics[height=100mm,width=165mm]{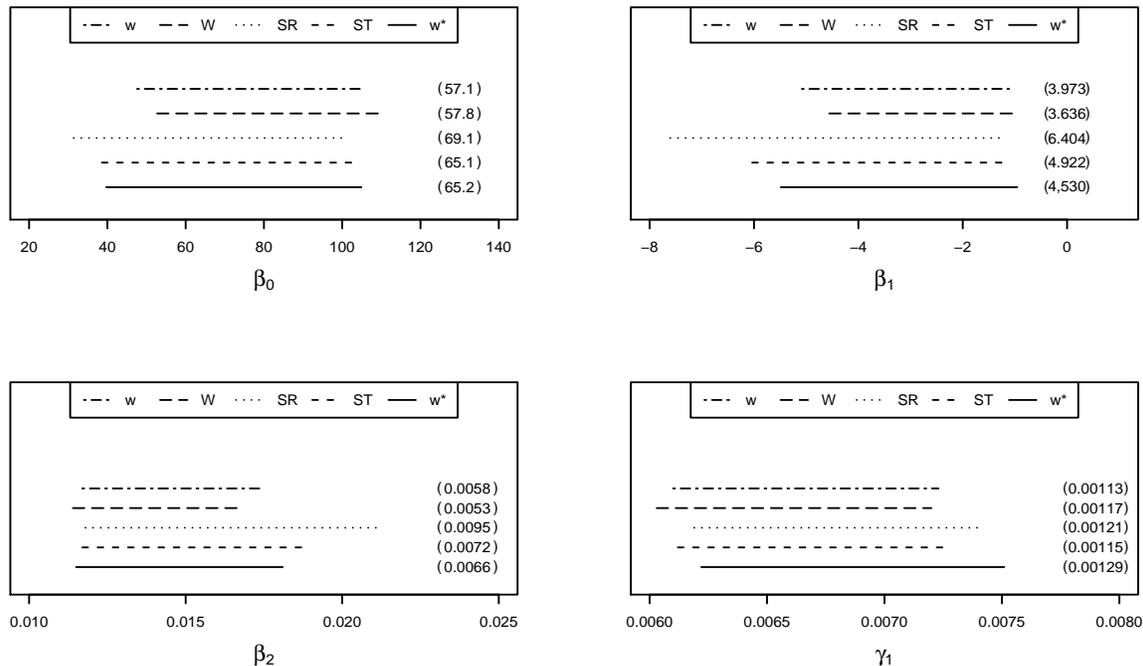}
  \caption{95\% confidence intervals; length in parentheses.}
  \label{fig:CI}
\end{figure}

The 95\% confidence intervals obtained by inverting the five tests are presented in Table~\ref{tab:CI} and in Figure~\ref{fig:CI}. It can be seen that the intervals obtained from the likelihood ratio and the Wald tests tend to be shorter than those obtained from the other tests as expected, since the first two are the most liberal. We
emphasize that theoretical and empirical findings indicate that the confidence intervals constructed from the adjusted likelihood ratio statistic should be regarded as the most accurate.
\begin{table}[!h]
\footnotesize
{
\begin{center}
\caption{95\% confidence intervals}\label{tab:CI}
\begin{tabular}
{l c c c c c}
\hline \hline
  &    $w$  &  $W$  & $S_R$ &  $S_T$  & $w^*$ \\
\hline 
$\beta_0$  &	(47.6; 104.7)	&	(52.7; 110.5)	&	(31.4; 100.5)	&	(38.6 ; 103.7)	&	(39.7; 104.9)	\\
$\beta_1$  &	($-5.084$; -1.111)	&	($-4.558$; $-0.922$)	&	($-7.604$; $-1.200$)	&	($-6.035$; $-1.113$)	&	($-5.486$; $-0.956$)	\\
$\beta_2$  &	(0.0117; 0.0175)	&	(0.0114; 0.0167)	&	(0.0118; 0.0213)	&	(0.0117; 0.0189)	&	(0.0115; 0.0181)	\\
$\gamma_1$ &	(0.00610; 0.00723)	&	(0.00603; 0.00720)	&	(0.00619; 0.00740)	&	(0.00612; 0.00727)	&	(0.00622; 0.00751)	\\
\hline \hline 
\end{tabular}
\end{center}
}
\end{table}
%\vspace{0.3cm}
%obtido com wheat-nãolinear-hetero-IC-beta0.ox, wheat-nãolinear-hetero-IC-beta1ox, wheat-nãolinear-hetero-IC-beta2.ox, wheat-nãolinear-hetero-IC-gama0.ox

Our second application deals with a dataset consisting of 34 men's decathlon performance at the 1988 Olympic Games
\citet[p. 304]{HAND}.\footnote{The dataset is
also available at {\tt http://www.stat.ncsu.edu/working\_groups/sas/sicl/data/olympic.dat}.}
%The final points scores of the competition for the 34 athletes are recorded in the ten events.
We assume that the score in high jump follows an extreme-value regression model (\ref{extreme-value}) 
with constant dispersion and systematic component for the location parameter given by
$$
\mu_t=\beta_0+\beta_1x_{t1}+\beta_2x_{t2}+\beta_3x_{t3}+\beta_4x_{t4}+\beta_5x_{t5},
$$
for $t=1,\ldots,34$. The covariates are the scores in the following events: javelin throw ($x_1$), 
long jump ($x_2$), discus throw ($x_3$), shot put ($x_4$), and pole vault ($x_5$).    
The statistics and the corresponding $p$-values for testing ${\mathcal H}_0: \beta_1=0$ against ${\mathcal H}_1: \beta_1 \neq 0$ 
are presented in Table~\ref{tab:decathlon_tests}.
\begin{table}[!h]
\begin{center}
\caption{Test statistics and the corresponding $p$-values}\label{tab:decathlon_tests}
\begin{tabular}
{c| c c c c c}
\hline \hline 
            & $w$  &  $W$  & $S_R$ &  $S_T$  & $w^*$\\
\hline        
statistic	& 4.0407 & 5.7161 & 2.8208  & 3.6293 &  2.6466\\
$p$-value     & 0.0444 & 0.0168 & 0.0930  & 0.0568 &  0.1038\\
\hline \hline 
\end{tabular}
\end{center}
\end{table}

Note that the $p$-values vary from $0.0168$ (Wald test) to $0.1038$ (adjusted likelihood ratio test). Only the likelihood ratio and Wald tests reject the null hypothesis at the 5\% nominal level, but this conclusion may be misleading since these tests
showed liberal behaviour in our simulations. The most reliable test, namely, the adjusted likelihood ratio test, does not suggest that
the null hypothesis should be rejected ($p$-value $= 0.1038$).

\section{Conclusion}\label{conclusion}

In this paper, we derived an adjusted version of the likelihood ratio statistic that
provides accurate inference in extreme-value regression models in small- to 
moderate-sized samples. Our simulation results suggest that the likelihood ratio and the Wald tests 
can be markedly oversized in small- and moderate-sized samples. 
The gradient test can be oversized, but much less than the other two tests.
The score test is even less size distorted and can be conservative in some cases.
We emphasize that the score test performs clearly better than the likelihood ratio, Wald and gradient tests, and is competitive with the adjusted likelihood ratio test in most cases. 
The adjusted likelihood ratio test obtained in this paper performs better than all the others. 
Although it requires more computational effort, is the least size distorted in 
most cases and it is, therefore, recommended for practical applications. 
We emphasize that our simulations were carried out in 
extreme-value regression models with linear and non-linear predictors for both location and dispersion parameters.
All simulation results exhibited reasonably similar behavior.

\section*{Acknowledgements}
We thank the referee for helpful comments and suggestions. We gratefully acknowledge the financial support 
from CNPq.

\small
\appendix
\section*{Appendix}
\label{ape:apendixqUpsilon}

Let $y \sim EV_{max}(\mu,\phi)$ and $z=(y-\mu)/\phi \sim EV_{max}(0,1)$. We have
$$
{\rm E}({z}^n \exp(-c{z}))
= \int_{-\infty}^{\infty}{z}^n \exp(-c{z})\exp(-{z}-\exp(-{z}))d{z}, \ \ n=0,1,\ldots.
$$
Using the transformation $v=\exp(-z)$, we have
\begin{eqnarray*}
{\rm E}({z}^n \exp(-c{z}))
&=& \int_{\infty}^{0}(-\ln v)^n v^c v \exp(-v)\frac{-1}{v}dv
= \int_{0}^{\infty}(-1)^n (\ln v)^n v^c \exp(-v)dv\\
&=& (-1)^n \Gamma ^{(n)}(1+c), \ \  n=1,2,\ldots; 
\end{eqnarray*}
see (4.358) in \cite{GRADSHTEYN}.
%$$
%\int_{0}^{\infty}x^{\nu-1}\exp(-\mu x)(\ln x)^n dx = \frac{\partial ^n}{\partial \nu^n}\{ \mu^{-\nu} \Gamma(\nu)\} \qquad {\rm para} \qquad n=0,1,2,3,\ldots
%$$
Since $\Gamma^{(1)}(n)=-(n-1)!(1/n+{\mathcal E}-\sum_{k=1}^n 1/k)$, we obtain
\begin{eqnarray*}
&&{\rm E}({z})= {\mathcal E}, \ \ {\rm E}({z}^2)=\Gamma^{(2)}(1),  \ \ {\rm E}(\exp(-c{z}))=\Gamma(1+c),\ \
{\rm E}(\exp(-(1+c){z}))=\Gamma(2+c), \\ 
&&{\rm E}({z} \exp(-{z}))={\mathcal E}-1, \ \ {\rm E}({z} \exp(-2{z}))=2{\mathcal E}-3, \ \ {\rm E}({z} \exp(-c{z}))=-\Gamma^{(1)}(1+c), \\ 
&&{\rm E}({z} \exp(-(1+c){z}))=-\Gamma^{(1)}(2+c), \ \ {\rm E}({z}^2 \exp(-{z}))=\Gamma^{(2)}(2),\\
&&{\rm E}({z}^2 \exp(-c{z}))=\Gamma^{(2)}(1+c), \ \ {\rm E}({z}^2 \exp(-(1+c){z}))=\Gamma^{(2)}(2+c).
\end{eqnarray*}
Let ${\mathnormal e}=(\exp(-(\phi_{11}/\phi_1){\mathfrak z}_1),\ldots,\exp(-(\phi_{1n}/\phi_n){\mathfrak z}_n))^{\!\top}$. Since 
$\Gamma^{(1)}(n)=\Gamma(n-1)+(n-1)\Gamma^{(1)}(n-1)$
and
$\Gamma^{(2)}(n)=2\Gamma^{(1)}(n-1)+(n-1)\Gamma^{(2)}(n-1)$ 
we can write after some calculations that
\begin{eqnarray*}
&&{\rm E}_{\omega }({\breve {\mathfrak z}})={\mathbf 1}, \ \
{\rm E}_{\omega }({\mathfrak z})={\mathcal E}{\mathbf 1}, \ \
{\rm E}_{\omega }({\mathcal Z})={\mathcal E}{\mathcal I}, \ \
{\rm E}_{\omega }({\breve {\mathcal Z}})={\mathcal I}, \ \
{\rm E}_{\omega }({\mathcal Z}{\breve {\mathfrak z}})=({\mathcal E}-1){\mathbf 1}, \\
&&{\rm E}_{\omega }({\mathcal Z}{\breve {\mathcal Z}})=({\mathcal E}-1){\mathcal I}, \ \
{\rm E}_{\omega }({\mathcal Z}^2{\breve {\mathcal Z}})=\Gamma^{(2)}(2){\mathcal I}, \ \
{\rm E}_{\omega}({\mathnormal e})= M{\mathbf 1},\ \
{\rm E}_{\omega }({\breve {\mathfrak z}}{\breve {\mathfrak z}}^{\!\top})={\mathbf 1} {\mathbf 1}^{\!\top} +{\mathcal I},\\
&&{\rm E}_{\omega }({\breve {\mathfrak z}}{\mathnormal e}^{\!\top})={\mathbf 1} {\mathbf 1}^{\!\top} M +CM,\ \
{\rm E}_{\omega }({\breve {\mathfrak z}}{\mathfrak z}^{\!\top})={\mathcal E}{\mathbf 1}{\mathbf 1}^{\!\top}-{\mathcal I},\ \
{\rm E}_{\omega }({\mathfrak z}{\breve {\mathfrak z}}^{\!\top})={\rm E}_{\omega }({\breve {\mathfrak z}}{\mathfrak z}^{\!\top}),\\
&&{\rm E}_{\omega }({\mathfrak z}{\mathnormal e}^{\!\top})={\mathcal E}{\mathbf 1}{\mathbf 1}^{\!\top}M-{\mathcal E}M-N,\ \
{\rm E}_{\omega }({\mathfrak z}{\mathfrak z}^{\!\top})={\mathcal E}^2{\mathbf 1} {\mathbf 1}^{\!\top} -{\mathcal E}^2{\rm I}+\Gamma^{(2)}(1){\mathcal I},\ \
{\rm E}_{\omega }({\mathcal Z}{\mathnormal e})=-N{\mathbf 1}, \\
&&{\rm E}_{\omega }({\breve {\mathfrak z}}{\mathnormal e}^{\!\top}{\mathcal Z}) =-{\mathbf 1} {\mathbf 1}^{\!\top}N - M - CN, \ \
{\rm E}_{\omega }({\mathfrak z}{\mathnormal e}^{\!\top}{\mathcal Z})=-{\mathcal E}{\mathbf 1} {\mathbf 1}^{\!\top}N + {\mathcal E}N + P, \\
&&{\rm E}_{\omega }({\mathcal Z}{\breve {\mathfrak z}}{\breve {\mathfrak z}}^{\!\top})=({\mathcal E}-1){\mathbf 1} {\mathbf 1}^{\!\top} +({\mathcal E}-2){\mathcal I}, \\
&&{\rm E}_{\omega }({\mathcal Z}{\breve {\mathfrak z}}{\mathnormal e}^{\!\top}) =({\mathcal E}-1){\mathbf 1}{\mathbf 1}^{\!\top} M - ({\mathcal E}-1) M - M -({\mathcal I}+C)N, \\
&&{\rm E}_{\omega }({\mathcal Z}{\breve {\mathfrak z}}{\mathfrak z}^{\!\top})={\mathcal E}({\mathcal E}-1){\mathbf 1} {\mathbf 1}^{\!\top} -{\mathcal E}({\mathcal E}-1){\mathcal I} + \Gamma^{(2)}(2){\mathcal I},\\
&&{\rm E}_{\omega }({\mathcal Z}{\breve {\mathfrak z}}{\mathnormal e}^{\!\top} {\mathcal Z})=-({\mathcal E}-1){\mathbf 1} {\mathbf 1}^{\!\top}N+({\mathcal E}-1)N + 2N + P + CP.
\end{eqnarray*}
Now, let $y_t \sim EV_{max}(\mu_{1t},\phi_{1t})$, ${\mathfrak z}_{1t}=(y_t-\mu_{1t})/\phi_{1t} \sim EV_{max}(0,1) $ and ${\mathfrak z}_t=(y_t-\mu_t)/\phi_t$. We can write ${\mathfrak z}_t=(\phi_{1t}/\phi_t) {\mathfrak z}_{1t}+(\mu_{1t}-\mu_t)/\phi_t$. Therefore,
\begin{equation}\label{z}
{\mathfrak z}=C{\mathfrak z}_1+D{\mathbf 1}
\end{equation}
and
\begin{equation}\label{e}
{\breve {\mathfrak z}}={\breve D}{\mathnormal e}_1,
\end{equation}
%\begin{equation}\label{Ze}
%{\mathcal Z}{\breve {\mathfrak z}}=C{\breve D}{\mathcal Z}_1{\mathnormal e}_1+D{\breve D}{\mathnormal e}_1,
%\end{equation}
where 
${\mathnormal e}_1=(\exp(-(\phi_{11}/\phi_1){\mathfrak z}_{11}),\ldots,\exp(-(\phi_{1n}/\phi_n){\mathfrak z}_{1n}))^{\!\top}$
%${\mathcal Z}_1={\rm diag}\{{\mathfrak z}_{11},\ldots,{\mathfrak z}_{1n}\}$
and the other quantities are as given earlier.

The vector $q$ given in (\ref{q}) is given by
$$
q=%{\rm E}_{\omega _1}[U(\omega _1)(l(\omega _1)-l(\omega ))]=
\left[
\begin{array}{c}
{\rm E}_{\omega _1}[U_\beta (\omega _1)\ell(\omega _1)]-{\rm E}_{\omega _1}[U_\beta (\omega_1)\ell(\omega )]\\
{\rm E}_{\omega _1}[U_\gamma (\omega _1)\ell(\omega _1)]-{\rm E}_{\omega _1}[U_\gamma (\omega _1)\ell(\omega _1)]
\end{array}
\right].
$$
From (\ref{loglikvet}), (\ref{scorebeta}), and the expected values obtained above, we have
\begin{eqnarray*}
{\rm E}_{\omega_1 }[U_\beta (\omega_1 )\ell (\omega_1 )]
&=&{\rm E}_{\omega_1 }\{X_1^{\!\top} \Phi_1^{-1} T_1 ({\mathbf 1}-{\breve {\mathfrak z}}_1)
[-{\mathnormal l}^{\!\top}-{\mathfrak z}_1^{\!\top}-{\breve {\mathfrak z}}_1^{\!\top}]{\mathbf 1}\}\nonumber\\
&=& X_1^{\!\top} \Phi_1^{-1} T_1\{-{\mathbf 1}{\mathnormal l}^{\!\top}-{\mathbf 1} {\rm E}_{\omega_1 }({\mathfrak z}_1^{\!\top})- {\mathbf 1}{\rm E}_{\omega_1 }({\breve {\mathfrak z}}_1^{\!\top})+ {\rm E}_{\omega_1 }({\breve {\mathfrak z}}_1){\mathnormal l}^{\!\top}+ {\rm E}_{\omega_1 }({\breve {\mathfrak z}}_1 {\mathfrak z}_1^{\!\top})+ {\rm E}_{\omega_1 }({\breve {\mathfrak z}}_1 {\breve {\mathfrak z}}_1^{\!\top}) \}{\mathbf 1}\nonumber\\
&=&X_1^{\!\top} \Phi_1^{-1} T_1\{-{\mathbf 1}\mathnormal {l}^{\!\top} -{\mathcal E}{\mathbf 1} {\mathbf 1}^{\!\top} -{\mathbf 1} {\mathbf 1}^{\!\top}+ {\mathbf 1} \mathnormal {l}^{\!\top} +({\mathcal E}{\mathbf 1} {\mathbf 1}^{\!\top}-{\mathcal I})+({\mathbf 1} {\mathbf 1}^{\!\top}+{\mathcal I})\}{\mathbf 1}=0{\mathbf 1}.
\end{eqnarray*}
Now, from (\ref{z}) and (\ref{e}), we have
\begin{eqnarray*}
{\rm E}_{\omega_1 }[U_\beta (\omega_1 )\ell (\omega )]
&=&{\rm E}_{\omega_1 }\{X_1^{\!\top} \Phi_1^{-1} T_1 ({\mathbf 1}-{\breve {\mathfrak z}}_1)
[-{\mathnormal l}^{\!\top}-{\mathfrak z}^{\!\top}-{\breve {\mathfrak z}}^{\!\top}]{\mathbf 1}\}\\
&=& X_1^{\!\top} \Phi_1^{-1} T_1{\rm E}_{\omega_1 }\{({\mathbf 1}-{\breve {\mathfrak z}}_1)
[-{\mathnormal l}^{\!\top}-{\mathfrak z}_1^{\!\top}C-{\mathbf 1}^{\!\top}D-{\mathnormal e}_1^{\!\top}{\breve D}]{\mathbf 1}\}\nonumber\\
&=&X_1^{\!\top} \Phi_1^{-1} T_1\{
-{\mathbf 1}\mathnormal {l}^{\!\top} 
-{\mathbf 1}{\rm E}_{\omega_1 }({\mathfrak z}_1^{\!\top})C 
-{\mathbf 1} {\mathbf 1}^{\!\top}D
- {\mathbf 1} {\rm E}_{\omega_1 }({\mathnormal e}_1^{\!\top}){\breve D} 
+{\rm E}_{\omega_1 }({\breve {\mathfrak z}}_1) \mathnormal {l}^{\!\top}
+{\rm E}_{\omega_1 }({\breve {\mathfrak z}}_1 {\mathfrak z}_1^{\!\top})C+\\
& &{\rm E}_{\omega_1 }({\breve {\mathfrak z}}_1){\mathbf 1}^{\!\top}D
+{\rm E}_{\omega_1 }({\breve {\mathfrak z}}_1{\mathnormal e}_1^{\!\top}){\breve D}    \}{\mathbf 1}\\
&=&X_1^{\!\top} \Phi_1^{-1} T_1
\{-{\mathbf 1}\mathnormal {l}^{\!\top} -{\mathcal E}{\mathbf 1} {\mathbf 1}^{\!\top}C -{\mathbf 1} {\mathbf 1}^{\!\top}D
-{\mathbf 1} ({\mathbf 1}^{\!\top}M ){\breve D} +{\mathbf 1}\mathnormal {l}^{\!\top}+({\mathcal E}{\mathbf 1}{\mathbf 1}^{\!\top}-{\mathcal I})C+\\
& &{\mathbf 1}{\mathbf 1}^{\!\top}D+
({\mathbf 1}{\mathbf 1}^{\!\top}M +CM){\breve D}\}{\mathbf 1}
=X_1^{\!\top} \Phi_1^{-1} T_1C(M{\breve D}-{\mathcal I}){\mathbf 1}.
\end{eqnarray*}
%onde $M$ é como definida no Capítulo~\ref{subcap:skovgaardmax}.
Hence,
$$
{\rm E}_{\omega _1}[U_\beta (\omega _1)l(\omega _1)]-{\rm E}_{\omega _1}[U_\beta (\omega _1)l(\omega )]
=X_1^{\!\top} \Phi_1^{-1} T_1C({\mathcal I}-M{\breve D}){\mathbf 1}.
$$

From (\ref{loglikvet}), (\ref{scoregamma}), and the results above,  we have
\begin{eqnarray*}
{\rm E}_{\omega_1 }[U_\gamma (\omega_1 )l(\omega_1 )]
&=&{\rm E}_{\omega_1 }\{Z_1^{\!\top} \Phi_1^{-1} H_1 ({\mathfrak z}_1-{\mathcal Z}_1{\breve {\mathfrak z}}_1-{\mathbf 1})[-\mathnormal {l}^{\!\top}-{\mathfrak z}_1^{\!\top}-{\breve {\mathfrak z}}_1^{\!\top}]{\mathbf 1}\}\nonumber\\
&=&Z_1^{\!\top} \Phi_1^{-1} H_1 
\{
-{\rm E}_{\omega_1 }({\mathfrak z}_1)\mathnormal {l}^{\!\top}
-{\rm E}_{\omega_1 }({\mathfrak z}_1{\mathfrak z}_1^{\!\top})
-{\rm E}_{\omega_1 }({\mathfrak z}_1{\breve {\mathfrak z}}_1^{\!\top})
+{\rm E}_{\omega_1 }({\mathcal Z}_1{\breve {\mathfrak z}}_1)\mathnormal {l}^{\!\top}\nonumber\\
&&+{\rm E}_{\omega_1 }({\mathcal Z}_1{\breve {\mathfrak z}}_1{\mathfrak z}_1^{\!\top})
+{\rm E}_{\omega_1 }({\mathcal Z}_1{\breve {\mathfrak z}}_1{\breve {\mathfrak z}}_1^{\!\top})
+{\mathbf 1} \mathnormal {l}^{\!\top}
+{\mathbf 1} {\rm E}_{\omega_1 }({\mathfrak z}_1^{\!\top})
+{\mathbf 1}{\rm E}_{\omega_1 }( {\breve {\mathfrak z}}_1^{\!\top})
\}{\mathbf 1} \\
&=&Z_1^{\!\top} \Phi_1^{-1} H_1 
\{
-\Gamma^{(2)}(1){\mathcal I}
+2{\mathcal E}{\mathcal I}
+\Gamma^{(2)}(2){\mathcal I}
-{\mathcal I}
\}{\mathbf 1} \\
&=&-Z_1^{\!\top} \Phi_1^{-1} H_1{\mathbf 1},
\end{eqnarray*}
and from (\ref{z}) and (\ref{e}) we get
\begin{eqnarray*}
{\rm E}_{\omega_1 }[U_\gamma (\omega_1 )l(\omega )]
&=&{\rm E}_{\omega_1 }\{Z_1^{\!\top} \Phi_1^{-1} H_1 ({\mathfrak z}_1-{\mathcal Z}_1{\breve {\mathfrak z}}_1-{\mathbf 1})[-\mathnormal {l}^{\!\top}-{\mathfrak z}^{\!\top}-{\breve {\mathfrak z}}^{\!\top}]{\mathbf 1}\}\nonumber\\
&=&{\rm E}_{\omega_1 }\{Z_1^{\!\top} \Phi_1^{-1} H_1 ({\mathfrak z}_1-{\mathcal Z}_1{\breve {\mathfrak z}}_1-{\mathbf 1})
[-\mathnormal {l}^{\!\top}-{\mathfrak z}_1^{\!\top}C-{\mathbf 1}^{\!\top}D-{\mathnormal e}_1^{\!\top}{\breve D}]{\mathbf 1}\}\nonumber\\
&=&
Z_1^{\!\top} \Phi_1^{-1} H_1 
\{
-{\rm E}_{\omega_1 }({\mathfrak z}_1){\mathnormal l}^{\!\top}
-{\rm E}_{\omega_1 }({\mathfrak z}_1{\mathfrak z}_1^{\!\top})C
-{\rm E}_{\omega_1 }({\mathfrak z}_1){\mathbf 1}^{\!\top}D
-{\rm E}_{\omega_1 }({\mathfrak z}_1{\mathnormal e}_1^{\!\top}){\breve D}\nonumber\\
&&
+{\rm E}_{\omega_1 }({\mathcal Z}_1{\breve {\mathfrak z}}_1)\mathnormal {l}^{\!\top}
+{\rm E}_{\omega_1 }({\mathcal Z}_1{\breve {\mathfrak z}}_1{\mathfrak z}_1^{\!\top})C
+{\rm E}_{\omega_1 }({\mathcal Z}_1{\breve {\mathfrak z}}_1){\mathbf 1}^{\!\top}D
+{\rm E}_{\omega_1 }({\mathcal Z}_1{\breve {\mathfrak z}}_1{\mathnormal e}_1^{\!\top}){\breve D}\nonumber\\
&&
+{\mathbf 1}\mathnormal {l}^{\!\top}
+{\mathbf 1}{\rm E}_{\omega_1 }({\mathfrak z}_1^{\!\top})C
+{\mathbf 1}{\mathbf 1}^{\!\top}D
+{\mathbf 1}{\rm E}_{\omega_1 }({\mathnormal e}_1^{\!\top}){\breve D}
\}{\mathbf 1}\nonumber\\
%&=&
%Z_1^{\!\top} \Phi_1^{-1} H_1 
%\{
%-{\mathcal E}{\mathbf 1}{\mathnormal l}^{\!\top}
%-({\mathcal E}^2{\mathbf 1} {\mathbf 1}^{\!\top} -{\mathcal E}^2{\rm I}+\Gamma^{(2)}(1){\mathcal I})C
%-{\mathcal E}{\mathbf 1}{\mathbf 1}^{\!\top}D\nonumber\\
%&&
%-({\mathcal E}{\mathbf 1}{\mathbf 1}^{\!\top}M-{\mathcal E}M-N){\breve D}
%+({\mathcal E}-1){\mathbf 1}\mathnormal {l}^{\!\top}\nonumber\\
%&&
%+({\mathcal E}({\mathcal E}-1){\mathbf 1} {\mathbf 1}{\!\top} -{\mathcal E}({\mathcal E}-1){\mathcal I} + \Gamma^{(2)}(2){\mathcal I})C
%+({\mathcal E}-1){\mathbf 1}{\mathbf 1}^{\!\top}D\nonumber\\
%&&
%+(({\mathcal E}-1){\mathbf 1}{\mathbf 1}^{\!\top}M - {\mathcal E}M -({\mathcal I}+C)N){\breve D}\nonumber\\
%&&
%+{\mathbf 1}\mathnormal {l}^{\!\top}
%+{\mathcal E}{\mathbf 1}{\mathbf 1}^{\!\top}C
%+{\mathbf 1}{\mathbf 1}^{\!\top}D
%+{\mathbf 1}(M{\mathbf 1})^{\!\top}{\breve D}
%\}{\mathbf 1}\nonumber\\
&=&
-Z_1^{\!\top} \Phi_1^{-1} H_1 C
\{{\mathcal E}{\mathcal I}+N{\breve D}\}{\mathbf 1}.
\end{eqnarray*}
It follows that
$$
{\rm E}_{\omega _1}[U_\gamma (\omega _1)l(\omega _1)]-{\rm E}_{\omega _1}[U_\gamma (\omega _1)l(\omega )]
%&=&-Z_1^{\!\top} \Phi_1^{-1} H_1{\mathbf 1}-(-Z_1^{\!\top} \Phi_1^{-1} H_1 C({\mathcal E}{\mathcal I}+N{\breve D} ){\mathbf 1} )\\
%&=&-Z_1^{\!\top} \Phi_1^{-1} H_1({\mathcal I}-C({\mathcal E}{\mathcal I}+N{\breve D} )){\mathbf 1}\\
=Z_1^{\!\top} \Phi_1^{-1} H_1(C({\mathcal E}{\mathcal I}+N{\breve D}) -{\mathcal I}){\mathbf 1}.
$$
Hence,
$$
q=
\left[
\begin{array}{c}
X_1^{\!\top} \Phi_1^{-1} T_1 C({\mathcal I}-M{\breve D} ){\mathbf 1} \\
Z_1^{\!\top} \Phi_1^{-1} H_1 \{C({\mathcal E}{\mathcal I}+N{\breve D})-{\mathcal I}\}{\mathbf 1}
\end{array}
\right].
$$

The matrix $\Upsilon$ given in (\ref{Upsilon}) can be written as
$$
\Upsilon =
\left[
\begin{array}{c c}
{\rm E}_{\omega _1}[U_\beta (\omega _1)U_\beta ^{\!\top}(\omega )] &
{\rm E}_{\omega _1}[U_\beta (\omega _1)U_\gamma (\omega )]\\
{\rm E}_{\omega _1}[U_\gamma (\omega _1)U_\beta ^{\!\top}(\omega )] &
{\rm E}_{\omega _1}[U_\gamma (\omega _1)U_\gamma ^{\!\top}(\omega )]
\end{array}
\right].
$$
From (\ref{scorebeta}), (\ref{e}), and the expected values obtained in the begining of this Appendix, we have
\begin{eqnarray*}
{\rm E}_{\omega _1}\{U_{\beta }(\omega _1)U_{\beta }^{\!\top}(\omega )\}
%&=&{\rm E}_{\omega _1}\{X_1^{\!\top} \Phi_1^{-1} T_1 ({\mathbf 1}-{\breve {\mathfrak z}}_1)
%[X^{\!\top} \Phi^{-1} T ({\mathbf 1}-{\breve {\mathfrak z}})]^{\!\top}\}\nonumber\\
&=&{\rm E}_{\omega _1}\{X_1^{\!\top} \Phi_1^{-1} T_1 ({\mathbf 1}-{\breve {\mathfrak z}}_1)
[X^{\!\top} \Phi^{-1} T ({\mathbf 1}-{\breve D}{\mathnormal e}_1)]^{\!\top}\}\nonumber\\
&=&X_1^{\!\top} \Phi_1^{-1} T_1 \bigl(
{\mathbf 1}{\mathbf 1}^{\!\top}
-{\mathbf 1}{\rm E}_{\omega _1}({\mathnormal e}_1^{\!\top}){\breve D}
-{\rm E}_{\omega _1}({\breve {\mathfrak z}}_1){\mathbf 1}^{\!\top}
+{\rm E}_{\omega _1}({\breve {\mathfrak z}}_1{\mathnormal e}_1^{\!\top}){\breve D}
\bigr)T \Phi^{-1}X\nonumber\\
&=&X_1^{\!\top} \Phi_1^{-1} T_1 
\bigl({\mathbf 1}{\mathbf 1}^{\!\top}-{\mathbf 1} {\mathbf 1}^{\!\top} M {\breve D}
-{\mathbf 1}{\mathbf 1}^{\!\top}
+({\mathbf 1} {\mathbf 1}^{\!\top} M +CM){\breve D}
\bigr)T \Phi^{-1}X\nonumber\\
&=&X_1^{\!\top} \Phi_1^{-1} T_1 CM{\breve D} T \Phi^{-1}X.
\end{eqnarray*}

The other blocks of $\Upsilon$ are derived in a similar fashion. We obtained
%
%From (\ref{scorebeta}), (\ref{scoregamma}) and the results above we have
%\begin{eqnarray*}
$$
{\rm E}_{\omega _1}\{U_{\beta }(\omega _1)U_{\gamma }^{\!\top}(\omega )\}
=X_1^{\!\top} \Phi_1^{-1} T_1 C\{{\mathcal I} + {\breve D}  (MD -M -CN)\}H\Phi^{-1} Z,
$$
%\end{eqnarray*}
%\begin{eqnarray*}
$$
{\rm E}_{\omega _1}\{U_{\gamma }(\omega _1)U_{\beta }^{\!\top}(\omega )\}
%&=&{\rm E}_{\omega _1}\{Z_1^{\!\top} \Phi_1^{-1} H_1 ({\mathfrak z}_1-{\mathcal Z}_1{\breve {\mathfrak z}}_1-{\mathbf 1})
%[X^{\!\top} \Phi^{-1} T ({\mathbf 1}-{\breve {\mathfrak z}})]^{\!\top} \}\nonumber\\
%&=&{\rm E}_{\omega _1}\{Z_1^{\!\top} \Phi_1^{-1} H_1 ({\mathfrak z}_1-{\mathcal Z}_1{\breve {\mathfrak z}}_1-{\mathbf 1})[X^{\!\top} \Phi^{-1} T ({\mathbf 1}-{\breve D}{\mathnormal e}_1)]^{\!\top} \}\nonumber\\
%&=&Z_1^{\!\top} \Phi_1^{-1} H_1 \{
%{\rm E}_{\omega _1}({\mathfrak z}_1){\mathbf 1}^{\!\top}
%-{\rm E}_{\omega _1}({\mathfrak z}_1{\mathnormal e}_1^{\!\top}){\breve D}
%-{\rm E}_{\omega _1}({\mathcal Z}_1{\breve {\mathfrak z}}_1){\mathbf 1}^{\!\top}
%+{\rm E}_{\omega _1}({\mathcal Z}_1{\breve {\mathfrak z}}_1{\mathnormal e}_1^{\!\top}){\breve D}
%-{\mathbf 1}{\mathbf 1}^{\!\top}\nonumber\\
%&&
%+{\mathbf 1}{\rm E}_{\omega _1}({\mathnormal e}_1^{\!\top}){\breve D}
%\}T \Phi^{-1} X \nonumber\\
%&=&Z_1^{\!\top} \Phi_1^{-1} H_1 \{
%{\mathcal E}{\mathbf 1}{\mathbf 1}^{\!\top}
%-({\mathcal E}{\mathbf 1}{\mathbf 1}^{\!\top}M-{\mathcal E}M-N){\breve D}
%-(({\mathcal E}-1){\mathbf 1}){\mathbf 1}^{\!\top}
%+(({\mathcal E}-1){\mathbf 1}{\mathbf 1}^{\!\top} M \nonumber\\
%&&
%- ({\mathcal E}-1) M - M -({\mathcal I}+C)N){\breve D}
%-{\mathbf 1}{\mathbf 1}^{\!\top}
%+{\mathbf 1}(M{\mathbf 1})^{\!\top}{\breve D}
%\}T \Phi^{-1} X \nonumber\\
=-Z_1^{\!\top} \Phi_1^{-1} H_1 CN{\breve D} T \Phi^{-1} X,
$$
%\end{eqnarray*}
%
%From (\ref{scoregamma}), (\ref{z}), (\ref{Ze}) and the results above  we have
%\begin{eqnarray*}
$${\rm E}_{\omega _1}\{U_{\gamma }(\omega _1)U_{\gamma }^{\!\top}(\omega )\}
=Z_1^{\!\top} \Phi_1^{-1} H_1 C\{{\mathcal E}{\mathcal I} + {\breve D} (N + CP - ND)\}H\Phi^{-1}Z. 
$$
%\end{eqnarray*}
%where $P$ é como definida no Capítulo~\ref{subcap:skovgaardmax} e lembrando que, ${\mathcal E}=-\Gamma^{(1)}(1)$.
Therefore,
$$
\Upsilon=
\left[
\begin{array}{c c}
X_1^{\!\top} \Phi_1^{-1} T_1 CM{\breve D} T \Phi^{-1}X&

X_1^{\!\top} \Phi_1^{-1} T_1 C\{{\mathcal I} + {\breve D}  (MD -M -CN)\}H\Phi^{-1} Z \\

-Z_1^{\!\top} \Phi_1^{-1} H_1 CN{\breve D} T \Phi^{-1} X  \quad&
Z_1^{\!\top} \Phi_1^{-1} H_1 C\{{\mathcal E}{\mathcal I} + {\breve D} (N + CP - ND)\}H\Phi^{-1}Z
\end{array}
\right].
$$
The vector $\overline q$ and the matrix $\overline \Upsilon$ are obtained from $q$ and $\Upsilon$ given above by replacing 
$X_1$, $\Phi_1$, $T_1$, $Z_1$, and $H_1$ by $\widehat X$, $\widehat \Phi$, $\widehat T$, $\widehat Z$, and $\widehat  H$, 
respectively, and
$X$, $\Phi$, $T$, $Z$, and $H$ by $\widetilde X$, $\widetilde \Phi$, $\widetilde T$, $\widetilde Z$, and $\widetilde  H$, 
respectively.

\end{document}